\documentclass{article}

\usepackage{amsmath,amsfonts,amsthm,amssymb,graphicx,tikz}
\usepackage[all,2cell,ps]{xy}

\theoremstyle{plain}
\newtheorem{thm}{Theorem}

\newtheorem{lem}[thm]{Lemma}
\newtheorem{prop}[thm]{Proposition}

\newtheorem{cor}[thm]{Corollary}

\newtheorem{proj}{Project}
\newtheorem{qtn}[proj]{Question}

\theoremstyle{definition}

\newtheorem{ex}[thm]{Example}

\newtheorem*{rem}{Remark}

\theoremstyle{remark}

\newcommand{\bbB}{\mathbb{B}}
\newcommand{\bbC}{\mathbb{C}}

\newcommand{\bbP}{\mathbb{P}}
\newcommand{\bbQ}{\mathbb{Q}}
\newcommand{\bbR}{\mathbb{R}}

\newcommand{\bbZ}{\mathbb{Z}}

\newcommand{\bfH}{\mathbf{H}}

\newcommand{\calC}{\mathcal{C}}
\newcommand{\calD}{\mathcal{D}}
\newcommand{\calE}{\mathcal{E}}

\newcommand{\calH}{\mathcal{H}}

\newcommand{\calO}{\mathcal{O}}

\newcommand{\calZ}{\mathcal{Z}}

\newcommand{\al}{\alpha}

\newcommand{\Gam}{\Gamma}
\newcommand{\del}{\delta}

\newcommand{\Del}{\Delta}

\newcommand{\Lam}{\Lambda}
\newcommand{\sig}{\sigma}
\newcommand{\Sig}{\Sigma}

\newcommand{\om}{\omega}

\DeclareMathOperator{\PSL}{PSL}
\DeclareMathOperator{\GL}{GL}
\DeclareMathOperator{\PGL}{PGL}

\DeclareMathOperator{\PU}{PU}
\DeclareMathOperator{\SU}{SU}

\DeclareMathOperator{\Gal}{Gal}

\newcommand{\ssm}{\smallsetminus}

\newcommand{\conj}{\overline}
\newcommand{\wh}{\widehat}
\newcommand{\wt}{\widetilde}

\newenvironment{pf}{\begin{proof}}{\end{proof}}

\newenvironment{enum}{\begin{enumerate}}{\end{enumerate}}

\usetikzlibrary{arrows}

\title{Cusp and $b_1$ growth for ball quotients and maps onto $\bbZ$ with finitely generated kernel}
\author{Matthew Stover\footnote{This material is based upon work supported by the National Science Foundation under Grant Number NSF 1361000. The author acknowledges support from U.S.\ National Science Foundation grants DMS 1107452, 1107263, 1107367 "RNMS: GEometric structures And Representation varieties" (the GEAR Network).} \\ \small{Temple University}\\ \small{\textsf{mstover@temple.edu}}}
\date{\today}

\begin{document}

\maketitle

\begin{abstract}
Let $M = \bbB^2 / \Gamma$ be a smooth ball quotient of finite volume with first betti number $b_1(M)$ and let $\calE(M) \ge 0$ be the number of cusps (i.e., topological ends) of $M$. We study the growth rates that are possible in towers of finite-sheeted coverings of $M$. In particular, $b_1$ and $\calE$ have little to do with one another, in contrast with the well-understood cases of hyperbolic $2$- and $3$-manifolds. We also discuss growth of $b_1$ for congruence arithmetic lattices acting on $\bbB^2$ and $\bbB^3$. Along the way, we provide an explicit example of a lattice in $\PU(2, 1)$ admitting a homomorphism onto $\bbZ$ with finitely generated kernel. Moreover, we show that any cocompact arithmetic lattice $\Gam \subset \PU(n, 1)$ of simplest type contains a finite index subgroup with this property.
\end{abstract}

\section{Introduction}\label{sec:Intro}

Let $\bbB^n$ be the unit ball in $\bbC^n$ with its Bergman metric and $\Gamma$ be a torsion-free group of isometries acting discretely with finite covolume. Then $M = \bbB^n / \Gamma$ is a manifold with a finite number $\calE(M) \ge 0$ of cusps. Let $b_1(M) = \dim H_1(M; \bbQ)$ be the first betti number of $M$. One purpose of this paper is to describe possible behavior of $\calE$ and $b_1$ in towers of finite-sheeted coverings. Our examples are closely related to the ball quotients constructed by Hirzebruch \cite{Hirzebruch} and Deligne--Mostow \cite{Deligne--Mostow, Mostow}.

For $n = 1$, $M$ is also known as either a (real) hyperbolic $2$-manifold or a hyperbolic quasiprojective curve, and the behavior of $b_1$ and $\calE$ is well-understood via the basic topology of Riemann surfaces. For $3$-manifolds, particularly finite volume hyperbolic $3$-manifolds, growth of $b_1$ in finite covering spaces has been a subject of immense interest \cite{CLR, CLR2, Venkataramana, Calegari--Dunfield, Agol}. For $2$- and $3$-manifolds that are the interior of a compact manifold, elementary algebraic topology completely describes the contribution of the cusps to $b_1$: for $2$-manifolds $b_1 \ge \calE - 1$, and for $3$-manifolds $b_1 \ge \calE$ (which is known by `half lives, half dies', and follows easily from Poincar\'e duality). The first result in this paper shows that no such theorem holds for two-dimensional ball quotients.

\begin{thm}\label{thm:MainCusped}
There exist infinite towers $\{A_j\}$, $\{B_j\}$, $\{C_j\}$, $\{D_j\}$ of distinct smooth finite volume noncompact quotients of $\bbB^2$ such that:
\begin{enum}

\item $b_1(A_j)$ is uniformly bounded for all $j$ and $\calE(A_j), \mathrm{vol}(A_j) \to \infty$ as $j \to \infty$;

\item $b_1(B_j)$ and $\calE(B_j)$ are both uniformly bounded but $\mathrm{vol}(B_j) \to \infty$ as $j \to \infty$;

\item $b_1(C_j)$, $\calE(C_j)$, and $\mathrm{vol}(C_j)$ all grow linearly with $j$;

\item there are positive constants $s, t$ and $0 < \sig < \tau < \upsilon < 1$ such that
\[
\mathrm{vol}(D_j)^{\sig} < s\, b_1(D_j) < \mathrm{vol}(D_j)^{\tau} < t\, \calE(D_j) < \mathrm{vol}(D_j)^{\upsilon}
\]
for all $j$.

\end{enum}
All these manifolds can be taken as quotients of $\bbB^2$ by an arithmetic lattice.
\end{thm}

This leads to the following question: \emph{For which nonnegative integers $\alpha, \beta$ does there exists a noncompact finite volume ball quotient with $b_1(M) = \alpha$ and $\calE(M) = \beta$?} It would be interesting to find infinitely many ball quotient manifolds with $b_1 = 0$ and arbitrarily many cusps, if such examples exist, especially in a tower of finite coverings; it seems that we do not know a single example for which the associated lattice is neat.

Our methods also require proving the following result, for which no example previously existed in the literature.

\begin{thm}\label{thm:fgkernel}
There exists a torsion-free lattice $\Gam$ in $\PU(2, 1)$ and a homomorphism $\rho : \Gam \to \bbZ$ with finitely generated kernel.
\end{thm}

This answers a question of Hersonsky--Paulin \cite{Hersonsky--Paulin} (see Theorem \ref{thm:MoreFG} below for a more definitive answer). Note that Theorem \ref{thm:fgkernel} is false for $\mathrm{PU}(1,1)$. Indeed, normal subgroups of finitely generated Fuchsian groups are either finite index or infinitely generated \cite{Greenberg}. Our first example is nonuniform. In joint work with Catanese, Keum, and Toledo \cite{CKST}, we will show that the fundamental group of the so-called Cartwright--Steger surface also has this property (in fact, we knew of that example first and the proof of Theorem \ref{thm:fgkernel} is modeled on the argument from \cite{CKST}). See \cite{DiCerbo--StoverComm} for more examples that build upon the one constructed in this paper.

Examples of lattices in $\PU(2, 1)$ admitting homomorphisms onto $\bbZ$ with infinitely generated kernel are known. Indeed, certain quotients of $\bbB^2$ admit a holomorphic fibration over a smooth projective curve of genus $g \ge 2$ (see \S \ref{sec:DM}), and this fibration induces homomorphisms onto $\bbZ$. Since $g \ge 2$, the kernel of the homomorphism from the fundamental group of the projective curve onto $\bbZ$ has infinitely generated kernel, and hence the associated normal subgroup of the lattice in $\PU(2,1)$ must also be infinitely generated. In fact, it follows from work of Napier and Ramachandran \cite{Napier--Ramachandran} that the kernel of a homomorphism to $\bbZ$ is infinitely generated if and only if, perhaps after passing to a finite-sheeted covering, there is an associated fibration over a hyperbolic curve. Studying this criterion more carefully, we will show the following.

\begin{thm}\label{thm:MoreFG}
Let $\Gam \subset \PU(n, 1)$ be a cocompact arithmetic lattice of `simplest type', i.e., one commensurable with $G(\bbZ)$, where $G$ is the algebraic group defined by a hermitian form on $\ell^{n + 1}$ for a CM field $\ell$. Then there exists $\Gam^\prime \subset \Gam$ of finite index and a surjection $\rho : \Gam^\prime \to \bbZ$ with finitely generated kernel.
\end{thm}

Passage to a subgroup of finite index in Theorem \ref{thm:MoreFG} can be a necessity. For example, for many small $n$ there are known arithmetic lattices in $\PU(n,1)$ generated by complex reflections, like those constructed by Mostow \cite{MostowReflect}, Deligne--Mostow/Mostow \cite{Deligne--Mostow, Mostow}, Allcock \cite{Allcock}, and Deraux--Parker--Paupert \cite{Deraux--Parker--Paupert}. All arithmetic groups generated by complex reflections arise from the construction considered in Theorem \ref{thm:MoreFG} \cite{StoverReflect}, and groups generated by complex reflections are generated by elements of finite order, so they have finite abelianization and hence no homomorphism onto $\bbZ$.

The proof of Theorem \ref{thm:MoreFG} uses a result of Clozel \cite{Clozelwedge} to ensure that $\bbB^n / \Gam^\prime$ is a manifold with holomorphic $1$-forms $\eta, \sig$ with $\eta \wedge \sig \neq 0$. This implies via \cite{Napier--Ramachandran} that the subspaces of $H^{1, 0}$ arising from homomorphisms to $\bbZ$ with finitely generated kernel all lie outside a finite union of proper subspaces. The existence of $\rho$ follows, though we cannot say much about finding a specific $\rho$. In contrast, Theorem \ref{thm:fgkernel} has a completely explicit proof that studies a fibration over an elliptic curve.

We also recall that for a (real) hyperbolic $3$-manifold $M$ and a surjection $\rho : \pi_1(M) \to \bbZ$, the kernel of $\rho$ is finitely generated if and only if $\rho$ is the homomorphism associated with a fibration of $M$ over the circle \cite{Stallings}. Fundamental groups of hyperbolic $3$-manifolds fibered over the circle and the lattices in Theorem \ref{thm:MoreFG} thus provide examples of hyperbolic groups admitting a homomorphism onto $\bbZ$ with finitely generated kernel. The lattices in this paper appear to be the first fundamental groups of a closed K\"ahler manifold with negative curvature with homomorphisms onto $\bbZ$ of this kind.

In a different direction, normal subgroups of lattices in $\bbR$-rank $1$ Lie groups have limit set the entire sphere at infinity $\partial X$, where $X$ is the associated symmetric space. Therefore, Theorem \ref{thm:MoreFG} has the following immediate consequence, for which appears to be new for $n \ge 4$.

\begin{cor}\label{cor:GIgroups}
For every $n \ge 2$, there exists a discrete, finitely generated, infinite covolume subgroup of $\PU(n, 1)$ with limit set the entire boundary $\partial \bbB^n$ of $\bbB^n$.
\end{cor}

Corollary \ref{cor:GIgroups} is false for $n = 1$, since all finitely generated Fuchsian groups with full limit set are lattices \cite[Thm.\ 4.5.1]{Katok}. For $n = 2,3$, the existence of examples were known as follows. Suppose $M$ is a smooth ball quotient fibered over a smooth projective curve $C$ of genus $g \ge 2$ (cf.\ \S \ref{sec:DM}). Then one can show that the kernel of the associated surjection $\pi_1(M) \to \pi_1(C)$ is finitely generated (see \cite{Kapovich}), and this kernel provides the desired examples. We also note that these groups are not finitely presented.

When $M = \bbB^2 / \Gamma$ is closed (that is, $\calE(M) = 0$) one can construct towers $\{M_j\}$ for which $b_1(M_j)$ is identically zero and other towers for which $b_1(M_j)$ grows linearly in the volume of $M_j$. Examples of towers with $b_1$ identically $0$ follow from work of Rogawski \cite{Rogawski}, and are quotients of the ball by a certain family of so-called \emph{congruence arithmetic} lattices (see \S \ref{ssec:Arithmetic}). For quotients of $\bbB^n$ for $n \ge 3$, a similar result follows from work of Clozel \cite{Clozel}. On the other hand, Marshall \cite{Marshall} showed that $b_1(M_j) \ll \mathrm{vol}(M_j)^{\frac{3}{8}}$ for all principal congruence arithmetic quotients of $\bbB^2$. Recall that $a_j \ll b_j$ for two nonnegative sequences $\{a_j\}, \{b_j\}$ when $a_j < c\, b_j$ for some constant $c > 0$ independent of $j$, and we say that $a_j \sim b_j$ when $a_j \ll b_j$ and $b_j \ll a_j$. Marshall then proved the following.

\begin{thm}[Marshall \cite{Marshall}]\label{thm:CongruenceGrowth}
There exists an arithmetic quotient $M_0$ of $\bbB^2$ and a tower of congruence coverings $\{M_j\}$ of $M_0$ such that
\[
\mathrm{vol}(M_j)^{\frac{3}{8}} \ll b_1(M_j).
\]
One can take $M_0$ to be compact or noncompact.
\end{thm}

In particular, the growth rate is exactly $\mathrm{vol}(M_j)^{\frac{3}{8}}$ in principal congruence towers, and one can take these manifolds to prove part (4) of Theorem \ref{thm:MainCusped}. The lattices used to prove Theorem \ref{thm:CongruenceGrowth} are precisely the lattices of `simple type', and the lattices to which Rogawski's vanishing results apply are those which are not of simple type. Marhsall's proof uses deep facts about the existence of certain automorphic representations. In the course of proving Theorem \ref{thm:MainCusped}, we will give a completely elementary and geometric proof of Theorem \ref{thm:CongruenceGrowth} by constructing towers of retractions of ball quotients onto totally geodesic complex curves. Our retractions are the ones described by Deraux \cite{Deraux}.

For $n = 3$, we can produce a similar result with exponent $\frac{3}{15}$. However our result is not optimal there, as Cossutta has a lower bound of $\mathrm{vol}(M_j)^{\frac{1}{4}}$ \cite{Cossutta}, and it would be interesting to achieve that lower bound geometrically. Our methods would easily imply that Theorem \ref{thm:CongruenceGrowth} holds for all $n$ and exponent $\frac{3}{n^2 + 2 n}$ if there exist examples of retractions of this kind for $n > 3$, but we know of no examples. Simon Marshall informed us that the correct upper bound from endoscopy should be $b_1(M_j) \ll \mathrm{vol}(M_j)^{\frac{n + 1}{n^2 + 2 n}}$ for all $n \ge 2$, and it would be very interesting to find the optimal growth rate and give a geometric interpretation.

We now describe the organization of the paper. In \S \ref{sec:Prelims} we briefly give some preliminaries on ball quotients and arithmetic groups. In \S \ref{sec:Hirzebruch} we describe Hirzebruch's ball quotient \cite{Hirzebruch} and build the families $\{A_j\}$ and $\{B_j\}$ from Theorem \ref{thm:MainCusped}. In the process, we prove Theorem \ref{thm:fgkernel}. We then prove Theorem \ref{thm:MoreFG} in \S \ref{sec:MoreFG}. In \S \ref{sec:DM} we construct the families $\{C_j\}$ and $\{D_j\}$ and prove Theorem \ref{thm:CongruenceGrowth}. In \S \ref{sec:Qns} we make some closing remarks and raise some questions.

\subsubsection*{Acknowledgments} We thank Simon Marshall and Domingo Toledo for several stimulating conversations related to various aspects of this paper. We also thank the referee whose comment led us to realize that we could prove Theorem \ref{thm:MoreFG}.

\section{Preliminaries}\label{sec:Prelims}

\subsection{The ball and its quotients}\label{ssec:Ball}

Let $\bbB^n$ denote the unit ball in $\bbC^n$ with its Bergman metric. See \cite{Goldman} for more on the geometry of $\bbB^n$. The holomorphic isometry group of $\bbB^n$ is $\PU(n, 1)$, and finite volume manifold quotients of $\bbB^n$ are quotients by torsion-free lattices in $\PU(n, 1)$. It is well-known that closed (resp.\ finite volume noncompact) manifold quotients of $\bbB^n$ are smooth projective (resp.\ quasiprojective) varieties. We also note that $\bbB^1$ is the Poincar\'e disk, and so the manifold quotients of $\bbB^1$ are precisely the finite volume (real) hyperbolic $2$-manifolds.

For what follows, let $\Gamma \subset \PU(n, 1)$ be a torsion-free lattice and $M = \bbB^n / \Gamma$. Since $\bbB^n$ has negative sectional curvature (in fact, constant holomorphic sectional curvature $-1$), $M$ has a finite number $\calE(M) \ge 0$ of topological ends. When $M$ is noncompact, each end of $M$ is homeomorphic to $N \times [0, \infty)$, where $N$ is an infranil $(n - 1)$-manifold (see e.g., \cite{McReynolds}).

The lattice $\Gamma$ is called \emph{neat} if the subgroup of $\bbC^*$ generated by its eigenvalues under the adjoint representation is torsion-free. In particular, a neat lattice is torsion-free. Every lattice contains a neat subgroup of finite index (see \cite[Ch.\ III]{AMRT}). Suppose that $\Gamma$ is a neat lattice and $M = \bbB^n / \Gamma$ is noncompact. Then $M$ admits a smooth toroidal compactification by $(n-1)$-dimensional abelian varieties generalizing the usual compactification of a punctured Riemann surface by a finite set of points \cite{AMRT, Mok}.

\subsection{Arithmetic lattices}\label{ssec:Arithmetic}

One method of constructing lattices in $\PU(n, 1)$ (in fact, the only method that is known to work for $n \ge 4$) is via arithmetic subgroups of algebraic groups. Let $G$ be a $\bbQ$-algebraic group with
\[
G(\bbR) \cong \SU(n, 1) \times \SU(n + 1)^r
\]
and $\Gam$ the image in $\PU(n, 1)$ of $G(\bbZ)$ after projection onto the $\SU(n, 1)$ factor of $G(\bbR)$ followed by projection of $\SU(n, 1)$ onto $\PU(n, 1)$. Then $\Gam$ is a lattice in $\PU(n,1)$. It is known that $\bbB^n / \Gam$ is noncompact if and only if $G$ is the special unitary group of a hermitian form of signature $(n, 1)$ over an imaginary quadratic field (e.g., see \cite[\S 3.1]{StoverVol}).

For any natural number $N$, we have the usual reduction homomorphisms
\[
r_N : G(\bbZ) \to G(\bbZ / N \bbZ)
\]
arising from inclusion of $G(\bbZ)$ into $\GL_m(\bbZ)$ for some $m$ and reducing matrix entries modulo $N$. Let $K_N$ be the kernel of $r_N$ and $\Gam(N)$ the image of $K_N$ in $\Gam = \Gam(1)$. We call $\Gam(N)$ a \emph{principal congruence kernel}. Any lattice $\Del$ in $\PU(n, 1)$ containing some $\Gam(N)$ as a subgroup of finite index is called a \emph{congruence arithmetic lattice}. We also note that the existence of neat subgroups of finite index is typically proven by showing that $\Gam(N)$ is neat for all sufficiently large $N$.

It is not hard to see that an arithmetic lattice in $\PU(n,1)$ with a homomorphism onto $\bbZ$ contains subgroups that are not congruence arithmetic \cite[p.\ 148]{Lubotzky--Segal}. Given an arbitrary arithmetic lattice $\Lam$, commensurable with $\Gam = \Gam(1)$ as above, we define congruence subgroups of $\Lam$ by
\[
\Lam(N) = \Lam \cap \Gam(N).
\]
In particular, we can talk of the family of congruence covers of any arithmetic quotient of $\bbB^n$.

\begin{ex}\label{ex:ArithmeticDef}
We recall the construction of the arithmetic lattices of $\PU(n, 1)$ of `simplest type'. In fact, all lattices appearing in this paper are of this kind. Let $k$ be a totally real field and $\ell$ a totally imaginary quadratic extension of $k$. We call $\ell$ a \emph{CM field}. We denote the nontrivial element of $\Gal(\ell / k)$ by $z \mapsto \conj{z}$, since this action extends to complex conjugation at any extension of a real embedding of $k$ to a complex embedding of $\ell$.

Let $h$ be a nondegenerate hermitian form on $\ell^{n + 1}$ with respect to the $\Gal(\ell / k)$-action. We obtain a $\bbQ$-algebraic group $G$ such that
\[
G(\bbQ) = \{g \in \PGL_{n + 1}(\ell)\ :\ {}^t \conj{g} h g = h \},
\]
where ${}^t \conj{g}$ is the $\Gal(\ell / k)$-conjugate transpose. For every real embedding $\nu : k \to \bbR$, $h$ extends to a hermitian form $h^\nu$ on $\bbC^{n + 1}$ whose signature is independent of the two complex embeddings of $\ell$ that extend $\nu$, and
\[
G(\bbR) \cong \prod_{\nu : k \to \bbR} \PU(h^\nu).
\]
We assume that there is a fixed $\nu_1$ such that $\PU(h^{\nu_1}) \cong \PU(n, 1)$ and that $\PU(h^\nu) \cong \PU(n + 1)$ for all $\nu \neq \nu_1$. Then projection onto the $\PU(n, 1)$ factor embeds $G(\bbZ)$ as a lattice which is cocompact if and only if $k \neq \bbQ$.
\end{ex}

\section{Hirzebruch's ball quotient}\label{sec:Hirzebruch}

The following construction is due to Hirzebruch \cite{Hirzebruch}. Let $\zeta = e^{\pi i / 3}$ and $\Lam = \bbZ[\zeta]$. Then $E = \bbC / \Lam$ is the elliptic curve of $j$-invariant $0$, and we let $S$ denote the abelian surface $E \times E$. Take coordinates $[z, w]$ on $S$ and consider the following elliptic curves on $S$:
\begin{align*}
C_0 &= \{w = 0\} \\
C_\infty &= \{z = 0\} \\
C_1 &= \{z = w\} \\
C_\zeta &= \{w = \zeta z\}
\end{align*}
Then $C_\al^2 = 0$ for each $\al$ and $C_\al \cap C_\beta = \{[0, 0]\}$ for every $\al \neq \beta$.

Consider the blowup $\wt{S}$ of $S$ at $[0,0]$ with exceptional divisor $D$, and let $\wt{C}_\al$ be the proper transform of $C_\al$ to $\wt{S}$. Then $\wt{C}_\al$ is an elliptic curve on $\wt{S}$ of self-intersection $-1$. If
\[
\calC = \bigcup_{\al \in \{0, \infty, 1, \zeta\}} \wt{C}_\al,
\]
define $M = \wt{S} \ssm \calC$. Hirzebruch conjectured and Holzapfel proved \cite{Holzapfel} that $M$ is the quotient of $\bbB^2$ by an arithmetic lattice $\Gam$. It is noncompact and $\wt{S}$ is the smooth toroidal compactification of the quotient of $\bbB^2$ by an arithmetic lattice. It appears as the third example in the appendix to \cite{StoverVol} and in \cite{Parker}. It is also one of the five complex hyperbolic manifolds of Euler number one that admits a smooth toroidal compactification \cite{DiCerbo--Stover}.

We need a few more facts before constructing the families $\{A_j\}$ and $\{B_j\}$ from Theorem \ref{thm:MainCusped}. Since the fundamental group is a birational invariant,
\[
\pi_1(\wt{S}) \cong \pi_1(S) \cong \bbZ^4.
\]
In particular, every \'etale cover $\pi : S^\prime \to S$ induces an \'etale cover $\wt{\pi} : \wt{S}^\prime \to \wt{S}$, and more specifically, $\wt{S}^\prime$ is the blowup of $S^\prime$ at the points over $[0, 0] \in S$. Every such cover is Galois (i.e., regular) with covering group $G$, a quotient of $\bbZ^4$.

Let $\wt{\pi} : \wt{S}^\prime \to \wt{S}$ be a finite \'etale cover with Galois group $G$ and degree $d$, so $\wt{S}^\prime$ is an abelian surface blown up at $d$ points. The inclusion $M \hookrightarrow \wt{S}$ induces a surjection
\[
\pi_1(M) \to \pi_1(\wt{S})
\]
\cite[Prop.\ 2.10]{Kollar}. The inverse image $N = \wt{\pi}^{-1}(M)$ of $M$ in $\wt{S}^\prime$ is an \'etale cover $M$, and the number of connected components of $N$ is the index of $\Gam = \pi_1(M)$ in $G$ under the induced homomorphism
\[
\pi_1(M) \to \pi_1(\wt{S}) \to G.
\]
This homomorphism is surjective, so $N$ is a connected covering of $M$. Lastly, if $\calD = \wt{\pi}^{-1}(\calC)$ is the pullback of $\calC$ to $\wt{S}^\prime$, then $N = \wt{S} \ssm \calD$ and $\calE(N)$ is the number of connected components of $\calD$. In other words, $N$ is the quotient of the ball by an arithmetic lattice for which the associated smooth toroidal compactification is an abelian surface blown up at $d$ points.

In order to compute $b_1$ for these coverings, we will use the following lemma due to Nori.

\begin{lem}[Nori \cite{NoriZar}]\label{lem:Nori}
Let $X$ and $Y$ be smooth connected varieties over $\bbC$ and $f : X \to Y$ an arbitrary morphism. Then:
\begin{itemize}

\item[A.] there is a nonempty Zariski-open $U \subset Y$ such that $f^{-1}(U) \to U$ is a fiber bundle in the usual topology.

\item[B.] if $f$ is dominant, the image of $\pi_1(X)$ has finite index in $\pi_1(Y)$.

\item[C.] if the general fiber $F$ of $f$ is connected and there is a codimension two subset $S$ of $Y$ outside which all the fibers of $f$ have at least one smooth point (i.e., $f^{-1}(p)$ is generically reduced on at least one irreducible component of $f^{-1}(p)$, for all $p \notin S$), then
\[
\pi_1(F) \to \pi_1(X) \to \pi_1(Y) \to 1
\]
is exact.

\end{itemize}
\end{lem}

We also use the following easy lemma from group theory.

\begin{lem}\label{lem:fgkernel}
Let $\Gam$ be a finitely generated group, $r \ge 1$, and $\rho : \Gam \to \bbZ^r$ a homomorphism with finitely generated kernel $K$. If $K$ can be generated by $n$ elements, then every finite index subgroup $\Gam^\prime \subseteq \Gam$ containing $K$ has abelianization of rank at most $n + r$.
\end{lem}

\begin{pf}
Let $\sig_1, \dots, \sig_n$ be generators for $K$. Then $\rho(\Gam^\prime) \subseteq \bbZ^r$ is a free abelian subgroup, and is hence generated by at most $r$ elements. Let $t_1, \dots, t_s$ ($s \le r$) be generators for $\rho(\Gam^\prime)$. Then we can lift $t_i$ to $\wh{t}_i \in \Gam^\prime$, and the $\wh{t}_i$ along with $\sig_1, \dots, \sig_n$ generate $\Gam^\prime$, hence $\Gam^\prime$ can be generated by $n + r$ elements. The cardinality of a generating set is an obvious upper bound for the rank of the abelianization of $\Gam^\prime$, so the lemma follows.
\end{pf}

Finally, we need some notation for $\pi_1(S)$. Write $\pi_1(S)$ as the group of translations of $\bbC^2$ generated by:
\begin{align*}
v_1 &= \begin{pmatrix} 1 \\ 0 \end{pmatrix} &v_3 = \begin{pmatrix} 0 \\ 1 \end{pmatrix} \\
v_2 &= \begin{pmatrix} \zeta \\ 0 \end{pmatrix} &v_4 = \begin{pmatrix} 0 \\ \zeta \end{pmatrix}
\end{align*}
Under the natural inclusions we then have:
\begin{align*}
\pi_1(C_0) &= \langle v_1, v_2 \rangle \\
\pi_1(C_\infty) &= \langle v_3, v_4 \rangle \\
\pi_1(C_1) &= \langle v_1 + v_3, v_2 + v_4 \rangle \\
\pi_1(C_\zeta) &= \langle v_1 + v_4, v_2 - v_3 + v_4 \rangle
\end{align*}
We now have all the tools necessary to build the families $\{A_j\}$ and $\{B_j\}$.

\subsubsection*{The towers $\{A_j\}$}

Fix a prime $p$ and consider the homomorphisms
\[
\rho_j : \pi_1(S) \cong \bbZ^4 \to P_j = \bbZ / p^j \bbZ
\]
such that $\rho_j(v_1)$ is a generator $\del_j$ for $P_j$ and $\rho_j(v_k)$ is trivial for $k \neq 1$. By the above discussion, this induces a connected $p^j$-fold cyclic cover $A_j$ of $M$. The restriction of $\rho_j$ to $\pi_1(C_\al)$ is surjective for $\al \in \{0, 1, \zeta\}$, so the cusp of $M$ associated with $\wh{C}_\al$ lifts to a single cusp of $A_j$. However, $\pi_1(C_\infty)$ is contained in the kernel of $\rho_j$, so the cusp of $M$ associated with $\wh{C}_\infty$ lifts to $p^j$ distinct cusps. Therefore, $\calE(A_j) = p^j + 3$.

We now must show that the betti numbers of the $A_j$ are uniformly bounded. Notice that each $\rho_j$ factors through the surjective homomorphism $\psi : \pi_1(M) \to \bbZ^2$ induced by the map $\wt{S} \to E$ given by blowdown followed projection of $S$ onto its first factor. By Lemma \ref{lem:fgkernel}, it suffices to show that $\psi$ has finitely generated kernel. In other words, we must prove the following.

\begin{prop}\label{prop:Hirzebruchfgkernel}
Let $E = \bbC / \bbZ[e^{2 \pi i / 3}]$, $S = E \times E$, $\wt{S}$ be the blowup of $S$ at $[0,0]$, and $M \subset \wt{S}$ be Hirzebruch's noncompact ball quotient. Consider the holomorphic map $M \to E$ induced by blowdown $\wt{S} \to S$ followed by projection of $S$ onto its first factor. Then the induced surjective homomorphism
\[
\psi : \pi_1(M) \to \pi_1(E) \cong \bbZ^2
\]
has finitely generated kernel.
\end{prop}

\begin{pf}
Notice that the general fiber $F$ of the map $M \to E$ is a reduced $3$-punctured elliptic curve. Indeed, $F$ intersects each $C_\al$ exactly once for $\al \in \{0, 1, \zeta\}$. The only fiber not of this type is the fiber above $0$, which is a $4$-punctured $\bbP^1$, namely the intersection of $M$ with the exceptional divisor $D$ of the blowup $\wt{S} \to S$. See Figure \ref{fig:FirstFibration}. This satisfies the assumptions of Lemma \ref{lem:Nori} with $X = M$ and $Y = E$, so there is an exact sequence
\[
\pi_1(F) \to \pi_1(M) \overset{\psi}{\to} \pi_1(E) \to 1.
\]
In particular, $\pi_1(F)$ maps onto $\ker(\psi)$. Since $\pi_1(F)$ is finitely generated, the proposition follows.
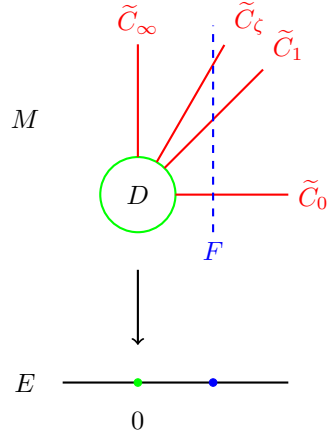
\begin{figure}
\begin{center}
\begin{tikzpicture}
\draw [thick, green] (0,0) circle (0.5);
\node at (0,0) {$D$};
\node at (-1.5, 1) {$M$};
\draw [thick, red, domain=0.5:2] plot (0, \x) node[above] {$\wt{C}_\infty$};
\draw [thick, red, domain=0.25:1.15] plot (\x, {sqrt(3)*\x}) node[above right] {$\wt{C}_\zeta$};
\draw [thick, red, domain=sqrt(2)/4:sqrt(2)+.25] plot (\x, \x) node[above right] {$\wt{C}_1$};
\draw [thick, red, domain=0.5:2] plot (\x, 0) node[right] {$\wt{C}_0$};
\draw [thick, <-] (0,-2) -- node[right]{} (0,-1);
\draw [thick, -] (-1,-2.5) -- (2,-2.5);
\node at (-1.5, -2.5) {$E$};
\draw [thick, dashed, blue] (1, -0.5) -- (1, 2.25);
\node [blue] at (1, -0.75) {$F$};
\draw [blue, fill=blue] (1, -2.5) circle (0.05cm);
\draw [green, fill=green] (0, -2.5) circle (0.05cm);
\node at (0, -3) {$0$};
\end{tikzpicture}
\caption{The fibration induced by projection onto the first factor of $S$.}\label{fig:FirstFibration}
\end{center}
\end{figure}
\end{pf}

Choosing the $\del_j = \rho_j(v_1)$ to be compatible with a family of homomorphisms $P_{j + 1} \to P_j$, the $\{A_j\}$ form a tower of coverings, which proves the first part of Theorem \ref{thm:MainCusped}, and the associated lattices in $\PU(2, 1)$ are arithmetic. The above also proves Theorem \ref{thm:fgkernel}.

\begin{pf}[Proof of Theorem \ref{thm:fgkernel}]
Combine Proposition \ref{prop:Hirzebruchfgkernel} with Lemma \ref{lem:fgkernel}.
\end{pf}

\begin{rem}
One can show directly that each $A_j$ has first betti number exactly $4$, and then avoid Lemmas \ref{lem:Nori} and \ref{lem:fgkernel}. Indeed, if $M$ is a smooth ball quotient admitting a smooth toroidal compactification $S$, then $H^1(M, \bbC) \cong H^1(S, \bbC)$ (see \cite{DiCerbo--StoverComm}). As $A_j$ has smooth toroidal compactification the blowup of an abelian surface, the betti number is exactly $4$. However, we chose to give the above proof because the argument is much more robust and may have applications in other situations where a direct betti number calculation is not possible for the smooth toroidal compactifications in the tower.
\end{rem}

\subsubsection*{The towers $\{B_j\}$}

Let $p$ be an odd prime, and fix a generator $\del_j$ for $P_j = \bbZ / p^j \bbZ$ (compatibly for a family of homomorphisms $P_{j + 1} \to P_j$). Since $p$ is odd, notice that $\del_j^2$ is also a generator for $P_j$. We then define homomorphisms $\rho_j : \pi_1(S) \to P_j$ by
\[
\rho_j(v_i) = \del_j, \quad i = 1,2,3,4.
\]
Let $\{B_j\}$ be the associated tower of $p^j$-fold cyclic covers of $M$. Here, the restriction of $\rho_j$ to $\pi_1(C_\al)$ is surjective for every $\al$, which implies each $\wt{C}_\al$ has a unique lift to the associated covering, and hence $\calE(B_j) = 4$.

Each $\rho_j$ factors through the homomorphism $\pi_1(S) \to \bbZ^2$ induced by the holomorphic mapping
\[
f([z, w]) = z + w
\]
from $S$ onto $E$. The general fiber of $f$ above $w \in E$ is the curve
\[
F_w = C_{-1} + [0, w]
\]
where $C_{-1}$ is the curve in coordinates $[z, w]$ on $S$ described by $\{w = -z\}$. Taking the proper transform $\wh{F}_w$ of $F_w$ to $\wt{S}$, we see that for $w \neq 0$, the induced map $M \to E$ has fiber above $w$ a punctured elliptic curve, i.e., $\wh{F}_w$ minus its points of intersection with each $\wh{C}_\al$, $\al \in \{0, \infty, 1, \zeta\}$. Above $0$, the fiber is the union of $D \cap M$ and $\wh{F}_0$ minus its intersection with each of the $\wh{C}_\al$. Each fiber is again reduced.

Lemma \ref{lem:Nori} again applies, so the kernel of the induced homomorphism from $\pi_1(M)$ to $\bbZ^2$ is finitely generated. Then $b_1(B_j)$ is uniformly bounded above by Lemma \ref{lem:fgkernel}. This proves the second part of Theorem \ref{thm:MainCusped}, and again the fundamental groups of our examples are arithmetic.

\section{The proof of Theorem \ref{thm:MoreFG}}\label{sec:MoreFG}

We begin with the following proposition.

\begin{prop}\label{prop:WedgeSubspaces}
Suppose that $V$ is a compact complex manifold that admits two holomorphic $1$-forms with nonzero wedge product. Given $\Lam \subset \PSL_2(\bbR)$ a cocompact lattice with $\calO = \bfH^2 / \Lam$ the associated closed hyperbolic $2$-orbifold and $f : V \to \calO$ a holomorphic mapping with connected fibers, define $W_f = f^*(H^{1, 0}(\calO)) \subset H^{1, 0}(V)$ to be the pull-backs to $V$ of holomorphic $1$-forms on $\calO$. Finally, define
\[
\calZ = \{ W_f \subset H^{1, 0}(M)\ :\ f : V \to \calO\ \textrm{as above} \}.
\]
Then $\calZ$ is a finite union of proper linear subspaces of $H^{1, 0}(V)$.
\end{prop}

\begin{pf}
It is clear that $\calZ$ is a union of linear subspaces that are maximal isotropic subspaces of $H^{1, 0}(V)$ for the wedge product. A version of the Castelnuovo--de Franchis theorem due to Catanese \cite[Thm.\ 1.10]{Catanese} implies that every maximal isotropic subspace $Z \in \calZ$ of dimension at least two determines a holomorphic mapping $f$ from $V$ onto a curve $C_Z$ of genus $g \ge 2$ and $Z = f^*(H^{1,0}(C))$. In particular, since $V$ supports two $1$-forms with nonzero wedge product, every such $Z \in \calZ$ is a proper subspace. The $Z_f \in \calZ$ associated with mappings onto hyperbolic $2$-orbifolds of genus $1$ are lines, and hence are also proper subspaces. The proposition follows immediately from the finiteness of the set of holomorphic maps $f$ from $V$ onto a hyperbolic $2$-orbifold \cite{Corlette--Simpson, Delzant}.
\end{pf}

We now prove Theorem \ref{thm:MoreFG}.

\begin{pf}[Proof of Theorem \ref{thm:MoreFG}]
Let $\Gam \subset \PU(n, 1)$ be a lattice satisfying the assumptions of the theorem. It is a theorem of Kazhdan that we can replace $\Gam$ with a torsion-free subgroup of finite index with infinite abelianization \cite[Thm.\ 2]{Kazhdan}. See also \cite[Thm.\ 15.2.1]{Bergeron--Clozel}. Then $\bbB^n / \Gam$ is an $n$-dimensional smooth complex projective variety (and a compact K\"ahler manifold) that supports a nontrivial holomorphic $1$-form $\eta$.

Using the action of Hecke operators, Clozel showed that there is $\Gam^\prime \subset \Gam$ of finite index and nontrivial holomorphic $1$-forms $\eta, \sig$ on $V = \bbB^n / \Gam^\prime$ such that $\eta \wedge \sig \neq 0$ \cite[Prop.\ 3.2]{Clozelwedge}. Indeed, since $\Gam$ has infinite abelianization, $\bbB^n / \Gam$ admits a holomorphic $1$-form $\om$. Our $1$-forms are then $\eta = (C * p^*\om)$ and $\sig = (C^\prime * p^* \om)$ in Clozel's notation, where $C, C^\prime$ are elements of the Hecke algebra $\calH_{\Gam^\prime}$ of $\Gam^\prime$, $p : \bbB^n / \Gam^\prime \to \bbB^n / \Gam$ is the covering map, and $\tau \mapsto C * \tau$ is the usual action of $C \in \calH_{\Gam^\prime}$ on $H^{1,0}(\bbB^n / \Gam^\prime)$. (Note that the cusp assumption is superfluous since our lattices are cocompact.) In particular, Proposition \ref{prop:WedgeSubspaces} applies to $V$.

Associated with every nonzero $\eta \in H^{1, 0}(V)$ there is a homomorphism
\[
\rho_\eta : \Gam^\prime \to \bbZ.
\]
If $\ker(\rho_\eta)$ is not finitely generated, it follows from work of Napier--Ramachandran that there exists a holomorphic map $f : V \to C$ onto a curve of genus $g \ge 1$ such that $\eta \in f^*(H^{1, 0}(C))$ \cite[Thm.\ 4.3]{Napier--Ramachandran}. Indeed, if $g \ge 2$, then $\eta$ clearly lies in the set $\calZ$ defined in the statement of Proposition \ref{prop:WedgeSubspaces}. However, this still holds when $C$ has genus $1$ since the induced homomorphism
\[
f_* :  \Gam^\prime \to \pi_1(C) \cong \bbZ^2
\]
factors through a surjective homomorphism $\Gam^\prime \to \Lam$, where $\Lam$ is a cocompact lattice in $\PSL_2(\bbR)$ (cf.\ \cite[Thm. 0.2]{Napier--Ramachandran2} and the remark following the statement), so again $\eta \in \calZ$.

In particular, if $\eta \notin \calZ$ then $\rho_\eta$ has finitely generated kernel. Since $\calZ$ is a finite union of proper subspaces of $H^{1, 0}(V)$, the existence of such an $\eta$ is immediate. This completes the proof.
\end{pf}

\begin{rem}
Note that the $1$-forms $\eta$ and $\sig$ with $\eta \wedge \sig \neq 0$ might not be a $1$-form associated with a finitely generated kernel. Indeed, both $\eta$ and $\sig$ might lie in distinct linear subspaces contained in $\calZ$.
\end{rem}

\begin{rem}
Using quasiprojective analogues of the above arguments, we expect that Theorem \ref{thm:MoreFG} also holds for nonuniform arithmetic lattices in $\PU(n, 1)$.
\end{rem}

\section{Deligne--Mostow orbifolds}\label{sec:DM}

To build the towers $\{C_j\}$ and $\{D_j\}$ and prove Theorem \ref{thm:CongruenceGrowth}, we use the ball quotient orbifolds constructed by Deligne and Mostow \cite{Deligne--Mostow, Mostow}. See \cite[\S 2]{Deraux} for an excellent introduction to the geometry and topology of these orbifolds; we use the notation of that paper in what follows. Given an $(n + 3)$-tuple $\mu = (\mu_i)$ of rational numbers, we consider the following condition:
\begin{equation}\label{eq:INT}
\big( 1 - \mu_i - \mu_j)^{-1} \in \bbZ \ \textrm{when}\ i \neq j\ \textrm{and}\ \mu_i + \mu_j < 1 \tag{INT}
\end{equation}
For any such $\mu$, Deligne and Mostow produced a finite volume ball quotient orbifold $\calO_\mu$ by a (sometimes partial) compactification of the space of $n + 3$ distinct points on $\bbP^1$. In other words, $\calO_\mu$ is the quotient of $\bbB^n$ by a lattice $\Lam_\mu$ that contains elements of finite order. Mostow then showed that one can relax \eqref{eq:INT} to a half-integral condition $\frac{1}{2}$\eqref{eq:INT} that produces new orbifolds. We refer to \cite{Deligne--Mostow, Mostow, KLW, Toledo, Deraux} for more on their geometry.

First, we need to recall when there are totally geodesic inclusions $\calO_\nu \hookrightarrow \calO_\mu$ between $m$- and $n$-dimensional Deligne--Mostow orbifolds. Let $\mu = (\mu_i)$ be an $(n + 3)$-tuple satisfying \eqref{eq:INT}, and (after reordering the $\mu_i$ if necessary), suppose that
\[
\sum_{j = m + 3}^{n + 3} \mu_j < 1.
\]
We then define the \emph{hyperbolic contraction}
\[
\nu = \big( \mu_1, \dots, \mu_{m + 2}, \sum_{j = m + 3}^{n + 3} \mu_j \big),
\]
which also satisfies \eqref{eq:INT}. Allowing the appropriate $n-m$ points to coincide then induces an inclusion of the space of $m + 3$ distinct points on $\bbP^1$ into the space of $n + 3$ distinct points, and this gives a totally geodesic embedding
\[
\calO_\nu = \bbB^m / \Lam_\nu \hookrightarrow \calO_\mu = \bbB^n / \Lam_\mu.
\]
There is a similar statement for $\frac{1}{2}$\eqref{eq:INT}. See \cite[\S 8]{Deligne--Mostow} or \cite{Deraux}.

Crucial to our construction are so-called \emph{forgetful maps}, which were studied in detail in \cite{Deraux}. Rather than inclusions of moduli spaces, we now consider when the surjective holomorphic map of moduli spaces given by forgetting some number of points induces a surjective holomorphic mapping between Deligne--Mostow orbifolds. Deraux \cite{Deraux} gave a complete treatment of how one determines the pairs $\mu, \nu$ for which this induces a surjective holomorphic mapping $\calO_\mu \to \calO_\nu$, and we refer the reader there for details.

What we exploit is that there are many instances where one can find a pair $(\mu, \nu)$ such that $\calO_\nu$ is both a geodesic suborbifold and a quotient space of $\calO_\mu$, and the associated composition
\[
\calO_\nu \hookrightarrow \calO_\mu \to \calO_\nu
\]
is the identity. In other words, $\calO_\mu$ admits a holomorphic \emph{retraction} onto $\calO_\nu$.

We will need to lift retractions to finite covering spaces, and must work in the category of (Riemannian) orbifolds and their orbifold fundamental groups (see \cite[III.1]{Bridson--Haefliger}). Let $X$ be a connected orbifold, $Y$ a suborbifold, and $f : X \to Y$ be a retraction of orbifolds. The induced map on (orbifold) fundamental groups $f_*$ then must satisfy
\[
f_*|_{\pi_1^{\mathrm{orb}}(Y)} = \mathrm{id}_{\pi_1^{\mathrm{orb}}(Y)}.
\]
We then have the following elementary lemma.

\begin{lem}\label{lem:RetractionGrowth}
Let $X$ be a connected orbifold, $Y$ be a suborbifold, and $f : X \to Y$ be a continuous retraction. Given a finite sheeted covering $X^\prime \to X$, let $Y^\prime$ be the finite covering of $Y$ associated with the finite index subgroup $f_*(\pi_1^{\mathrm{orb}}(X^\prime))$ of $\pi_1^{\mathrm{orb}}(Y)$. Then $Y^\prime$ lifts to a submanifold of $X^\prime$, and $X^\prime$ admits a retraction onto $Y^\prime$.
\end{lem}

\begin{pf}
Let $h : X^\prime \to Y$ be projection $X^\prime \to X$ followed by $f$. Then we have a commutative diagram
\[
\begin{matrix}
Y^\prime & \hookrightarrow & X^\prime & \to & Y^! \\
\downarrow & & \downarrow & & \downarrow & \\
Y & \hookrightarrow & X & \to & Y
\end{matrix}
\]
where $Y^!$ is the minimal regular cover of $Y$ through which $h$ factors. Then
\[
\deg(Y^! \to Y) = [\pi_1^{\mathrm{orb}}(Y) : f_*(\pi_1^{\mathrm{orb}}(X^\prime))].
\]
We want to show that $Y^! = Y^\prime$. Notice that $Y^\prime$ must map onto $Y^!$, since $h(Y^\prime) = Y$, so
\[
\deg(Y^! \to Y) \le \deg(h|_{Y^\prime}).
\]

We now prove the opposite inequality. The restriction of $h$ to $Y^\prime$ is the covering $Y^\prime \to Y$ followed by the identity map $Y \to Y$, hence
\[
\deg(h|_{Y^\prime}) = \deg(f|_Y) = [\pi_1^{\mathrm{orb}}(Y) : \pi_1^{\mathrm{orb}}(Y^\prime)] = [f_*(\pi_1^{\mathrm{orb}}(Y)) : f_*(\pi_1^{\mathrm{orb}}(Y^\prime))],
\]
as $f_*$ is the identity on $\pi_1^{\mathrm{orb}}(Y)$. Since $\pi_1^{\mathrm{orb}}(Y^\prime)$ is a subgroup of $\pi_1^{\mathrm{orb}}(X^\prime)$ we certainly have
\[
[\pi_1^{\mathrm{orb}}(Y) : f_*(\pi_1^{\mathrm{orb}}(Y^\prime))] = [\pi_1^{\mathrm{orb}}(Y) : \pi_1^{\mathrm{orb}}(Y^\prime)] \le [\pi_1^{\mathrm{orb}}(Y) : f_*(\pi_1^{\mathrm{orb}}(X^\prime))].
\]
This all combines to give
\[
\deg(h|_{Y^\prime}) \le \deg(Y^! \to Y),
\]
and the lemma follows.
\end{pf}

\begin{rem}
We also could have stated the lemma in terms of retractions onto totally geodesic subspaces of locally symmetric spaces. One then replaces $\pi_1^{\mathrm{orb}}$ with the appropriate lattice in a Lie group and proceeds with exactly the same argument. However, in either situation it does not suffice to merely pass to a manifold cover of $X$ and argue with manifolds, as one needs the lemma to know that the manifold cover retracts onto the appropriate geodesic subspace!
\end{rem}

It follows immediately from the lemma that $b_1(X^\prime) \ge b_1(Y^\prime)$. We now explain how one constructs the families $\{C_j\}$ and $\{D_j\}$.

\subsubsection*{The towers $\{C_j\}$}

Let $\calO_\mu = \bbB^2 / \Gam_\mu$ be a cusped Deligne--Mostow orbifold and $f : \calO_\mu \to \calO_\nu$ be a surjective holomorphic mapping onto $\calO_\nu = \bbB^1 / \Gam_\nu$. For example, one can take
\begin{align*}
\mu &= \left( \frac{2}{6}, \frac{2}{6}, \frac{3}{6}, \frac{4}{6}, \frac{1}{6} \right) \\
\nu &= \left( \frac{1}{6}, \frac{3}{6}, \frac{4}{6}, \frac{4}{6} \right)
\end{align*}
While $\calO_\mu$ has cusps, we note that $\calO_\nu$ is compact.

Then, we can find a neat subgroups $\Gam \subset \Gam_\mu$ and $\Del \subset \Gam_\nu$ of finite index for which we have a surjective holomorphic mapping
\[
h : C_0 = \bbB^2 / \Gam \to \Sig = \bbB^1 / \Del
\]
for which the induced map $h_* : \Gam \to \Del$ is also surjective. Indeed, first replace $\Gam_\nu$ with a neat subgroup $\Del_0$ of finite index and take $\Gam_0 \subset \Gam_\mu$ to be the inverse image of $\Del_0$ in $\Gam_\mu$ under the surjection $\Gam_\mu \to \Gam_\nu$. We then take $\Gam$ to be a neat subgroup of finite index in $\Gam_0$ and set $\Del$ to be image of $\Gam$ in $\Del_0$, which is neat since $\Del_0$ is neat. Note that $\Sig$ is a compact Riemann surface of genus $g \ge 2$.

Fix a cusp of $C_0$ and let $P \subset \Gam$ be a representative for the conjugacy class of maximal parabolic subgroups of $\Gam$ associated with this cusp. We claim that $h_*(P)$ is a (possibly trivial) cyclic subgroup of $\Sig$. Indeed, $P$ is nilpotent, so $h_*(P)$ is also nilpotent. However, $\Sig$ is a hyperbolic surface group, so every nilpotent subgroup is cyclic.

Given a finite covering $\Sig^\prime \to \Sig$, Galois with group $G$, recall that the number of cusps of the induced covering $C^\prime$ of $C$ over our chosen cusp of $C$ is exactly the index of $P$ in $G$ under the associated surjection $\Gam \to G$. Since $b_1(\Sig) \ge 4$ and $h_*(P)$ is cyclic, we can find a homomorphism $\rho : \Del \to \bbZ$ with $h_*(P) \subset \mathrm{ker}(\rho)$. Let $\rho_j$ be the composition of $\rho$ with reduction modulo $j$ and let $\Del_j$ be the kernel of $\rho_j$.

Set:
\begin{align*}
\Sig_j &= \bbB^1 / \Del_j \\
C_j &= \bbB^2 / h_*^{-1}(\Del_j)
\end{align*}
Then $\Sig_j$ (resp.\ $C_j$) is a covering of $\Sig$ (resp.\ $C_0$) of degree $j$. Furthermore, $b_1(\Sig_j)$ grows linearly in $j$, hence $b_1(C_j)$ does as well. Since $P \subset h_*^{-1}(\Del_j)$, we see that the cusp associated with $P$ lifts to $j$ cusps of $C_j$. In particular, $\calE(C_j)$ also grows linearly in $j$. Therefore $\{C_j\}$ has the required properties.

\begin{rem}
Taking
\begin{align*}
\mu &= \left( \frac{2}{6}, \frac{2}{6}, \frac{3}{6}, \frac{3}{6}, \frac{1}{6}, \frac{1}{6} \right) \\
\nu &= \left( \frac{1}{6}, \frac{3}{6}, \frac{4}{6}, \frac{4}{6} \right)
\end{align*}
one also sees that there exist quotients of $\bbB^3$ for which Theorem \ref{thm:MainCusped}(3) holds. That is, the $3$-dimensional complex orbifold $\calO_\mu$ retracts onto the $1$-dimensional complex orbifold $\calO_\nu$, and we can thus construct a tower of coverings of $\calO_\mu$ for which the first betti number and number of cusps grows linearly in the covering degree.
\end{rem}

\subsubsection*{Theorem \ref{thm:CongruenceGrowth} and the towers $\{D_j\}$}

Let $D_0 = C_0 = \bbB^2 / \Gam$ and $\Sig = \bbB^1 / \Del$ be as above, and notice that our map is, in fact, a retraction onto a geodesic submanifold. We now let $D_j = \bbB^2 / \Gam(N_j)$ be the congruence covering of $D_0$ of level $N_j$ as in \S \ref{ssec:Arithmetic}, where $\{N_j\}$ is an infinite sequence of integers to be determined soon. Then
\[
\Del(N_j) = \Del \cap \Gam(N_j)
\]
defines a family of congruence subgroups of $\Del$. It follows from Lemma \ref{lem:RetractionGrowth} that $b_1(D_j) \ge b_1(\Sig_j)$, where $\Sig_j$ is now the congruence covering $\bbB^1 / \Del(j)$ of $\Sig$.

Since $\# \PSL_2(\bbZ / N_j \bbZ) \sim N_j^3$, using strong approximation for Zariski-dense subgroups \cite{Weisfeiler, Nori}, one can find an appropriate infinite sequence $\{N_j\}$ such that $b_1(\Sig_j)$ grows like $N_j^3$, hence $b_1(D_j)$ grows at least like $N_j^3$ in this tower. Indeed, strong approximation implies that reduction modulo $N_j$ is surjective for an infinite sequence of integers $N_j$. However, $G(\bbZ / N_j \bbZ)$ has order $N_j^8$, so a similar strong approximation argument shows that
\[
\mathrm{vol}(D_j) \sim [\Gam : \Gam(N_j)] \sim \#G(\bbZ / N_j \bbZ) \sim N_j^8,
\]
and we see that
\[
\mathrm{vol}(D_j)^{\frac{3}{8}} \ll b_1(D_j).
\]
This reproduces the lower bound from \cite{Marshall}. If we take $D_j$ to be the quotient of $\bbB^2$ by a congruence arithmetic lattice, an upper bound with the same exponent follows from \cite{Marshall}.

To calculate $\calE(D_j)$, we need to compute the index in $G(\bbZ / N_j \bbZ)$ of a given maximal parabolic subgroup of $\Gam$. If $P$ is a maximal parabolic subgroup of $\Gam$ then its image in $\Gam / \Gam(N_j)$ is contained in some Borel subgroup, but since the diagonal entries of parabolic matrices in $\Gam$ are units of a fixed imaginary quadratic number field, we see that their images in $G(\bbZ / N_j \bbZ)$ are, up to a fixed multiplicative factor, contained in the subgroup of strictly upper-triangular matrices. Elementary counting implies that this index is of order $N_j^5$. As argued above for $b_1$, it follows that $\calE(D_j) \sim \mathrm{vol}(D_j)^{5/8}$.

This completes the proof for noncompact quotients when $n = 2$. Indeed, we can take
\begin{align*}
\sig &= \frac{3}{8} - \epsilon \\
\tau &= \frac{3}{8} + \epsilon \\
\upsilon &= \frac{5}{8} + \epsilon
\end{align*}
for any $0 < \epsilon < \frac{1}{4}$.

Applying the same argument to
\begin{align*}
\mu &= \left( \frac{2}{6}, \frac{2}{6}, \frac{3}{6}, \frac{3}{6}, \frac{1}{6}, \frac{1}{6} \right) \\
\nu &= \left( \frac{1}{6}, \frac{3}{6}, \frac{4}{6}, \frac{4}{6} \right)
\end{align*}
gives noncompact examples when $n = 3$, where now $[\Gam : \Gam(N_j)]$ grows like $N_j^{15}$, so $b_1$ grows at least as fast as $N_j^{1/5}$. We note, however, that $\calE$ now grows like $N_j^{2/5}$. For compact examples, consider:
\begin{align*}
\mu_1 &= \left( \frac{3}{8}, \frac{3}{8}, \frac{3}{8}, \frac{7}{8} \right) \\
\mu_2 &= \left( \frac{3}{8}, \frac{3}{8}, \frac{3}{8}, \frac{3}{8}, \frac{4}{8} \right) \\
\mu_3 &= \left( \frac{1}{8}, \frac{3}{8}, \frac{3}{8}, \frac{3}{8}, \frac{3}{8}, \frac{3}{8} \right)
\end{align*}
By \cite[Thm.\ 3.1(v)]{Deraux}, there is a holomorphic retraction $f : \calO_{\mu_3} \to \calO_{\mu_1}$ induced by a forgetful map. Clearly the restriction of $f$ to $\calO_{\mu_2}$ is a retraction of it onto $\calO_{\mu_1}$. These give all the desired examples.

\begin{rem}
Notice that our lattice $\Gam$ is not necessarily congruence arithmetic, since the initial Deligne--Mostow lattice $\Gam_\mu$ might not be congruence arithmetic. One can often show that the quotient $\calO_{\Sig \mu}$ of $\calO_\mu$ by a certain symmetric group is the quotient of the ball by a maximal arithmetic lattice, which is then congruence arithmetic. Our noncompact example in dimension $2$ used to construct the tower $\{C_j\}$ also has the property that the map $\calO_\mu \to \calO_\nu$ descends to $\calO_{\Sig \mu}$. For $n = 3$, we do not know that the families $\{M_j\}$ in Theorem \ref{thm:CongruenceGrowth} are necessarily congruence arithmetic, though we suspect that such examples exist.
\end{rem}

\begin{rem}
Unfortunately, this result is not optimal for $n = 3$, as Cossetta proved the lower bound $\mathrm{vol}(M_j)^{\frac{1}{4}}$ for the growth of $b_1$ \cite{Cossutta}. It would be interesting to see if a more subtle use of fibrations over curves could meet, or beat, that bound.
\end{rem}

\section{Closing questions and remarks}\label{sec:Qns}

We begin reiterating a question from the introduction.

\begin{qtn}
Fix $n \ge 2$. For which pairs $(\al, \beta)$ is there a smooth finite volume quotient $M = \bbB^n / \Gam$ with $b_1(M) = \al$ and $\calE(M) = \beta$? Can we always assume $\Gam$ is arithmetic and/or neat?
\end{qtn}

Particularly interesting is the case $(0, \beta)$ for any $\beta \ge 1$. We do not know of an infinite family of examples, though \cite{StoverVol} contains examples for which the lattice is arithmetic but not neat. We apparently do not know a single example where the lattice is neat.

\begin{qtn}
Does there exist a manifold quotient $M$ of $\bbB^n$ with $\calE(M) = 1$? For which $n$ can $M$ be the quotient by a neat and/or arithmetic lattice?
\end{qtn}

No examples are known for $n \ge 2$. In an earlier paper \cite{StoverCusps}, we showed that for any $k \ge 1$, there exists a constant $n_k$ such that $\calE(M) > k$ for every arithmetic quotient of $\bbB^n$ with $n \ge n_k$. In particular, for $n$ sufficiently large, if one-cusped quotients of $\bbB^n$ exist, they cannot be arithmetic.

\begin{qtn}
Let $M$ be a finite volume quotient of $\bbB^n$ and $\{M_j\}$ a tower of finite-sheeted coverings. What are the possible growth rates for $b_1(M_j)$?
\end{qtn}

As mentioned in the introduction, we know examples for all $n \ge 2$ where $b_1(M_j)$ is identically zero \cite{Rogawski, Clozel}. There are also examples where the growth is nontrivial \cite{Kazhdan, Venkataramana}. In both families, the lattices can be chosen to be arithmetic.

\begin{qtn}
Let $M$ be the quotient of $\bbB^n$ by a congruence arithmetic lattice, $n \ge 2$, and $\{M_j\}$ a family of congruence coverings. What are the possible growth types for $b_1(M_j)$ as a function of the covering degree?
\end{qtn}

When $M$ is noncompact, one can always determine the growth rate of $\calE(M_j)$ by elementary counting methods. Simon Marshall informed us that endoscopy should give a bound $b_1(M_j) \ll \mathrm{vol}(M_j)^{\frac{n + 1}{n^2 + 2 n}}$ for principal congruence lattices. Cossutta has the upper and lower bounds $\frac{n+2}{(n+1)^2}$ and $\frac{n - 2}{(n + 1)^2}$, respectively. If $M$ retracts onto a holomorphically embedded geodesic submanifold $\bbB^1 / \Lam$, then the methods in this paper give towers with a lower bound of the form $\frac{3}{n^2 + 2n}$ (and could do more if the number of distinct retractions increases with $j$), but we do not know a single example of a quotient of $\bbB^n$, $n \ge 4$, that retracts onto a geodesic suborbifold of any codimension.

\begin{qtn}
Let $M$ be a quotient of $\bbB^n$, and suppose that $M$ contains a totally geodesic quotient of $\bbB^m$ for some $1 \le m < n$. Does $M$ admit a finite sheeted covering $N \to M$ such that $N$ retracts onto one of its $m$-dimensional geodesic suborbifolds?
\end{qtn}

Long and Reid showed that, for any $3 \le m \le n$, a complex hyperbolic $n$ manifold cannot retract onto a totally geodesic \emph{real} hyperbolic $m$-submanifold \cite{Long--ReidRetract}. The analogous result is always true for hyperbolic $3$-manifolds and totally geodesic hyperbolic $2$-manifolds \cite{Agol}. There are also a number of known examples of hyperbolic $n$-manifolds retracting onto totally geodesic hyperbolic $m$-submanifolds \cite{Bergeron--Haglund--Wise}, though the general case for $n \ge 4$ remains open.

\bibliography{B1growth}

\end{document}